\documentclass[reqno]{amsart}

\usepackage{amsmath}
\usepackage{amssymb}
\usepackage{amsthm}
\usepackage{xspace}
\usepackage{pictexwd}
\usepackage{float}  \restylefloat{figure} 
\usepackage{enumerate}

\theoremstyle{plain}
\newtheorem{thm}{Theorem}
\newtheorem*{conj}{Conjecture}
\newtheorem{prop}[thm]{Proposition}
\newtheorem{lem}[thm]{Lemma}
\newtheorem{cor}[thm]{Corollary}

\theoremstyle{remark}
\newtheorem*{rem}{Remark}

\def\forcehmode{\hskip0pt\relax}

\let\myskip=\medskip

\def\definebb#1=#2.{\def#1{{{\mathbb #2}^{\vphantom{x}}}}}
\definebb \c=C.
\definebb \r=R.
\definebb \z=Z.
\definebb \q=Q.

\def\wA{w_A}
\def\wB{w_B}

\def\tinf{t_{\infty}}
\def\cinf{c_{\infty}}
\def\tA{t_A}
\def\tAn{t_A^{\circ}}
\def\cA{c_A}

\def\tBp{t_B^{+}}
\def\tBm{t_B^{-}}
\def\tBpm{t_B^{\pm}}

\def\chp{H_{\c}^2}
\def\rhp{H_{\r}^2}
\def\cdo{\c^{2,1}}

\def\lll{$[\ell_1,\ell_2,\ell_3]$}
\def\ppp{$(p_1,p_2,p_3)$}
\def\Gappp{\Ga(p_1,p_2,p_3)}
\def\Gpppn{G(p_1,p_2,p_3;n)}

\def\ppi{$(p,p,\infty)$}

\def\iii{$(\infty,\infty,\infty)$}

\def\uui{$(4,4,\infty)$}
\def\Guuin{G(4,4,\infty;n)}

\def\uuu{$(4,4,4)$}
\def\Guuun{G(4,4,4;n)}

\def\dd{\partial}
\def\bz{\bar z}
\def\moddpi{~\mod2\pi}

\def\al{{\alpha}}
\def\be{{\beta}}
\def\Ga{{\Gamma}}
\def\ga{{\gamma}}
\def\de{{\delta}}

\def\la{{\lambda}}
\def\si{{\sigma}}

\let\emptyset=\varnothing
\def\st{\,\,\big|\,\,}
\def\<{\langle}
\def\>{\rangle}
\def\ie{i.e.\xspace}

\def\hcross{\boxtimes}

\let\ge=\geqslant
\let\le=\leqslant

\def\defit{\it}

\DeclareMathOperator{\ctan}{ctan}

\DeclareMathOperator{\id}{id}

\DeclareMathOperator{\PU}{PU}
\DeclareMathOperator{\SU}{SU}
\DeclareMathOperator{\trace}{trace}
\DeclareMathOperator{\U}{U}

\let\Re=\undefined \DeclareMathOperator{\Re}{Re}
\let\mod=\undefined \DeclareMathOperator{\mod}{mod}

\begin{document}

\author[Anna Pratoussevitch]{Anna Pratoussevitch}
\address{Mathematisches Institut\\ Universit\"at Bonn\\ Beringstra{\ss}e~1 \\ 53115 Bonn}
\email{anna@math.uni-bonn.de}

\title{Traces in Complex Hyperbolic Triangle Groups}

\begin{date}  {\today} \end{date}

\thanks{Research partially supported by NSF grant DMS-0072607 
and by SFB 611 of the DFG}

\begin{abstract} 
We present several formulas for the traces of elements
in complex hyperbolic triangle groups generated by complex reflections.

The space of such groups of fixed signature is of real dimension one.
We parameterise this space by a real invariant~$\alpha$ of triangles
in the complex hyperbolic plane.
The main result of the paper is a formula, which expresses the trace
of an element of the group as a Laurent polynomial in $e^{i\alpha}$
with coefficients independent of $\alpha$ and computable using a certain
combinatorial winding number.
We also give a recursion formula for these Laurent polynomials
and generalise the trace formulas for the groups generated
by complex $\mu$-reflections.

We apply these formulas to prove some discreteness and some non-dis\-cre\-te\-ness
results for complex hyperbolic triangle groups.

 \end{abstract}










\maketitle

\section{Introduction}

We study representations of real hyperbolic triangle groups,
\ie groups generated by reflections in the sides of triangles in~$\rhp$,
in the holomorphic isometry group $\PU(2,1)$ of the complex hyperbolic plane~$\chp$.

\myskip
We use the following terminology:
A {\defit complex hyperbolic triangle\/} is a triple $(C_1,C_2,C_3)$ of complex geodesics in~$\chp$.
If the complex geodesics $C_{k-1}$ and $C_{k+1}$ meet at the angle $\pi/p_k$
we call the triangle $(C_1,C_2,C_3)$ a {\defit \ppp-triangle}.
If the complex geodesics $C_{k-1}$ and $C_{k+1}$ are ultra-parallel
(\ie their closures in $\chp\cup\dd\chp$ do not intersect) with distance $\ell_k$
we call the triangle $(C_1,C_2,C_3)$ a {\defit \lll-triangle}.
A {\defit complex geodesic\/} (or {\defit complex slice\/}) in~$\chp$ is the fixed point set of a complex reflection.
A {\defit complex reflection\/} is an element of $\PU(2,1)$ conjugate to the map
$$[z_1:z_2:z_3]\mapsto[z_1:-z_2:-z_3].$$
For more details on complex reflections and complex and real slices see section~\ref{cpxhypplane}.

\myskip
We call a subgroup of $\PU(2,1)$ generated by complex reflections $\iota_k$ in the sides~$C_k$
of a complex hyperbolic \ppp-triangle $(C_1,C_2,C_3)$ a {\defit \ppp-triangle group\/}.
We call a representation of the group
$$\Gappp=\<\ga_1,\ga_2,\ga_3\st\ga_k^2=(\ga_{k-1}\ga_{k+1})^{p_k}=1~\hbox{for all}~k\in\{1,2,3\}\>,$$
where $\ga_{k+3}=\ga_k$, and the relation $(\ga_{k-1}\ga_{k+1})^{p_k}=1$ is to omit for $p_k=\infty$,
into the group $\PU(2,1)$ given by taking the generators $\ga_k$ of $\Gappp$
to the generators $\iota_k$ of a \ppp-triangle group a {\defit \ppp-representation\/}.
In the ultra-parallel case we define similarly \lll-{\defit triangle groups\/} and \lll-{\defit representations\/}. 

\myskip
For fixed \ppp\ resp. $(\ell_1,\ell_2,\ell_3)$ the space of complex hyperbolic triangle groups is of real dimension one.
There is a canonical path $\rho_t$ of representations.
The starting point $\rho_0$ for the path is the case when the triangle lies in a real slice.
Real slices are the second type of totally geodesic subspaces in~$\chp$.

\myskip
For \ppp-triangle groups, according to Richard Schwartz~\cite{Sch2002}, the conjectural picture is as follows:

\begin{conj}
We assume $p_1\le p_2\le p_3$.
Define $\wA=\iota_3\iota_2\iota_3\iota_1$ and $\wB=\iota_1\iota_2\iota_3$.

\myskip
A \ppp-representation is a discrete embedding if and only if neither $\wA$ nor $\wB$ is elliptic.
The set of the corresponding parameter values is a closed symmetric interval.

\myskip
If the element $\wA$ becomes elliptic before $\wB$,
we say that the triple \ppp\ is of {\defit type A},
else we say that the triple is of {\defit type B}.

\myskip
The triple \ppp\ is of type A if $p_1<10$ and of type B if $p_1>13$.

\myskip
If the triple \ppp\ is of type A then there is a countable collection of parameters,
for which the \ppp-representation is infinite and discrete but not injective.
If the triple \ppp\ is of type B then there are no such discrete but not
injective \ppp-representations.
\end{conj}

R.~Schwartz proved this conjecture for \iii-groups in~\cite{Sch2001:dented},
and for \ppp-groups with $p_1$, $p_2$, $p_3$ sufficiently large in~\cite{Sch2003:dehn}.
For sufficiently large $p_1$, $p_2$, $p_3$ the triple \ppp\ is of type B.

\myskip
The case of \ppi-triangle groups was studied in~\cite{W} by Justin Wyss-Gallifent.
It turns out that the triple \ppi\ is of type~A for $p\le13$ and of type~B for $p\ge14$.
He also described some of the discrete but not injective \uui-representations.

\myskip
In general, to prove that a certain \ppp-representation is discrete and
injective, it is necessary and sufficient to show that the images of elements
of $\Gappp$ of infinite order are not regular elliptic.
Here the traces come in, namely, an element of $\SU(2,1)$ is not regular
elliptic if and only if the value of certain discriminant function (introduced
by W.~Goldman) on the trace of the element is not negative.

\myskip
The main results of this paper are formulas for traces of elements in
complex hyperbolic triangle groups, a combinatorial trace formula
(theorem~\ref{comb}) and a recursive trace formula (theorem~\ref{reccur})
as well as their applications.
These formulas generalise the results of Hanna Sandler~\cite{S95}
on ideal triangle groups.

\myskip
We parameterise the one-dimensional space of \ppp-triangle groups by an invariant~$\alpha$
of triangles in the complex hyperbolic plane.
The combinatorial trace formula express the trace of an element of the group
as a Laurent polynomial in $e^{i\alpha}$ with coefficients independent
of $\alpha$ and computable using a certain combinatorial winding number.
Then we give a recursion formula for these Laurent polynomials.
We also generalise the trace formulas for the groups generated
by complex $\mu$-reflections.

\myskip
We apply the formulas to prove some discreteness and some non-discreteness results
for complex hyperbolic triangle groups.

\myskip
For instance, we compute the parameter value $\tA$ such that for $|t|>\tA$ the
corresponding \ppp-representations are not discrete embeddings because the
element~$\wA$ is regular elliptic (proposition~\ref{def-tA} and corollary~\ref{cor-tA}).
Furthermore, for a certain subfamily of triangle groups we can partially
confirm the conjecture of Schwartz about type~A and ~B.
We find a sufficient condition for a triple \ppp\ to be of type~B.
This condition implies for example that triples \ppi\ with $p\ge14$ and
triples $(p,2p,2p)$ with $p\ge12$ are of type~B (proposition~\ref{typeB}).
We also obtain similar results for ultra-parallel triangle groups.

\myskip
The paper is organised as follows:
In section~\ref{cpxhypplane} we recall the basic notions of complex hyperbolic
geometry, specially for the complex hyperbolic plane~$\chp$.
In sections~\ref{cpxhyptri} and~\ref{anglinvcompar} we define angular
invariant~$\al$, which classifies complex hyperbolic triangles with fixed
angles up to isometry, and we compare this invariant with invariants defined
for some special cases by U.~Brehm, J.~Hakim and H.~Sandler.
Section~\ref{winding} contains the description of certain combinatorial
functions of words, first of all the winding number.

\myskip
After that we are prepared to state and to prove in
sections~\ref{combtraceformula} and ~\ref{rectraceformula} our main results,
the combinatorial trace formula (thm.~\ref{comb}) and the recursive trace formula (thm.~\ref{reccur}). 
In section~\ref{examples} we discuss as the first application of the recursive
trace formula the resulting formulas for the traces of short words in triangle groups.
We prove in section~\ref{someproperties} some properties of the traces of
elements in a triangle group, which are also used later in section~\ref{appl-disc}.
In section~\ref{mu-refl} we describe the generalisations of our trace formulas 
for groups generated by $\mu$-reflections (for the definition see section~\ref{cpxhypplane}).

\myskip
Sections~\ref{appl-non-disc}--\ref{appl-disc} contain applications of the
trace formulas.
In section~\ref{appl-non-disc} we discuss necessary conditions for a triangle
group representation to be a discrete embedding.
Section~\ref{spec-family} deals with a subfamily of \ppp-triangle groups with
the property $\frac{\pi}{p_1}=\frac{\pi}{p_2}+\frac{\pi}{p_3}$
resp.\ of \lll-triangle groups with the property $\ell_3=\ell_1+\ell_2$.
The triangle groups in this subfamily seem to share a lot of properties of the
ideal triangle groups.
For this subfamily we can prove more precise statements about type~A and~B.
In section~\ref{appl-disc} we discuss some arithmetic properties of the traces for certain (finite) family of signatures
and parameter values of triangle groups.
Finally, in section~\ref{hist-geogr} we summarise some of the known results on
complex hyperbolic triangle groups and describe them in terms of our parameter
$r_1,r_2,r_3$, and $\al$ (for notation compare section~\ref{cpxhyptri}).

\myskip
Part of this work was done during an enjoyable stay at the University of Maryland.
I am grateful to William Goldman and Richard Schwartz for their kind invitation
and for helpful discussions
and to the University of Maryland for its hospitality.
I would like to thank
Werner Ballmann,
Egbert Brieskorn,
Martin Deraux,
Ilya Dogolazky,
Tadeusz Januszkiewicz,
John Parker,
Blake Pelzer,
and Gregor Weingart
for useful conversations related to this work.

\section{Complex Hyperbolic Plane}

\label{cpxhypplane}

In this section we recall some basic notions of complex hyperbolic geometry.
The general references on complex hyperbolic geometry are \cite{G99}, 
\cite{P03}, and \cite{BH99}.

\myskip
{\bf Complex Hyperbolic Plane:}
Let $\cdo$ denote the vector space $\c^3$ equipped with the Hermitian form
$$\<z,w\>=z_1\bar w_1+z_2\bar w_2-z_3\bar w_3$$
of signature $(2,1)$. We call a vector $z\in\cdo$ {\it negative}, {\it
null\/}, or {\it positive}, according as $\<z,z\>$ is negative, zero, or positive.
Let $P(\cdo)$ denote the projectivisation of~$\cdo-\{0\}$.
We denote the image of $z\in\cdo$ under the projectivisation map by $[z]$.
We also write $[z]=[z_1:z_2:z_3]$ for $z=(z_1,z_2,z_3)\in\cdo$.
The {\it complex hyperbolic plane\/} $\chp$ is the projectivisation of the set of negative
vectors in $\cdo$. Its ideal boundary $\dd\chp$ is defined as the
projectivisation of the set of null vectors in $\cdo$.
The complex hyperbolic plane $\chp$ is a K\"ahler manifold
of constant holomorphic sectional curvature.
The holomorphic isometry group of $\chp$ is the projectivisation $\PU(2,1)$
of the group $\SU(2,1)$ of complex linear transformations, which preserve the Hermitian form.

\myskip
The {\bf Hermitian cross product:} $\hcross:\cdo\times\cdo\to\cdo$ is defined by
$$
  z\hcross w
  =\begin{pmatrix} -1&0&0\\ 0&-1&0\\ 0&0&1 \end{pmatrix}(\bar z\times\bar w)
  =\begin{pmatrix}\,\overline{z_3w_2-z_2w_3}\,\\\,\overline{z_1w_3-z_3w_1}\,\\ \,\overline{z_1w_2-z_2w_1}\,\end{pmatrix}.
$$
The Hermitian cross product of two vectors yields a vector perpendicular to both of them
$$\<z\hcross w,z\>=\<z\hcross w,w\>=0.$$
Furthermore we have
$$\<a\hcross c,b\hcross c\>=\overline{\<a,c\>\<c,b\>-\<a,b\>\<c,c\>},$$
in particular
$$\<a\hcross b,a\hcross b\>=|\<a,b\>|^2-\<a,a\>\<b,b\>.$$

\myskip
{\bf Totally Geodesic Submanifolds:}
There are two kinds of totally geodesic submanifolds of dimension~$2$ in
$\chp$, {\it complex slices} (or {\it complex geodesics}) and {\it real
slices} (or {\it totally real totally geodesics subspaces}).
Complex geodesics are obtained by projectivisation of $2$-dimensional complex
subspaces of $\cdo$.
Given any two points in $\chp$, there is a unique complex geodesic
containing them.
Any positive vector $c\in\cdo$ determines
a $2$-dimensional complex subspace 
$$\{z\in\cdo\st\<c,z\>=0\}$$
and a complex geodesic, which is the projectivisation of this subspace.
The vector $c$ is called a {\it polar vector} of the complex geodesic.
A polar vector can be normalised to $\<c,c\>=1$.
Conversely, any complex geodesic is represented by a polar vector.
A typical example is the complex slice $\{[z:0:1]\in\chp\}$ with polar vector $c=(0,1,0)^T$.
Any complex slice is isometric to this one.

\myskip
We now describe possible configurations of two complex geodesics (compare proposition~6.8 in~\cite{P03}).
For two complex geodesics one of the following situations occur:
\begin{enumerate}[$\bullet$]
\item they coincide,
\item they intersect in a single point in~$\chp$,
\item their closures in $\chp\cup\dd\chp$ intersect in a single point in $\dd\chp$,
\item their closures in $\chp\cup\dd\chp$ are disjoint.
\end{enumerate}

\myskip
Let $C_1$ and $C_2$ be two complex geodesics with normalised polar vectors $c_1$ and $c_2$ respectively.
The complex geodesics $C_1$ and $C_2$ intersect in $\chp$ if and only if $|\<c_1,c_2\>|<1$,
in which case $|\<c_1,c_2\>|=\cos\angle(C_1,C_2)$,
where $\angle(C_1,C_2)$ is the angle of intersection between $C_1$ and $C_2$,
and $c_1\hcross c_2$ is a negative vector, corresponding to the intersection point.
The angle of intersection can be defined as
$$\angle(C_1,C_2)=\min\{\angle(\la_1c_1,\la_2c_2)\st\la_1,\la_2\in\r\},$$
where the angles on the right hand side are measured in the underlying real vector space of $\cdo$.
The angle satisfies $\angle(C_1,C_2)\in[0,\pi/2]$.
The complex geodesics $C_1$ and $C_2$ intersect in $\dd\chp$ if and only if $|\<c_1,c_2\>|=1$,
in this case $c_1\hcross c_2$ is a null vector, corresponding to the intersection point.
The closures of the complex geodesics $C_1$ and $C_2$ are disjoint if and only
if $|\<c_1,c_2\>|>1$, in which case $|\<c_1,c_2\>|=\cosh(\ell/2)$,
where $\ell$ is the distance between $C_1$ and $C_2$.

\myskip
Any real slice is isometric to $\{[z:w:1]\in\chp\st z,w\in\r\}$.
Real slices are fixed point sets of {\it real reflections,}
\ie antiholomorphic isometries conjugate to the map
$$[z:w:1]\mapsto[\bar z:\bar w:1].$$

\myskip
{\bf Complex Reflections:}
Given a complex geodesic~$C$, there is a unique isometry~$\iota_C$ in $\PU(2,1)$
of order~$2$, whose fixed point set equals~$C$.
We call this isometry the {\it complex reflection\/} in~$C$
(or the {\it inversion\/} on~$C$).
The complex reflection in $C$ is represented by an element $\iota_C\in\SU(2,1)$
that is given by
$$\iota_C(z)=-z+2\frac{\<z,c\>}{\<c,c\>}c,$$
where $c$ is a polar vector of~$C$.

\myskip
{\bf Complex $\mu$-Reflections:}
For a unit complex number $\mu$ a complex $\mu$-reflection is an element of
$\PU(2,1)$ conjugate to the map
$$[z_1:z_2:z_3]\mapsto[\mu z_1:z_2:z_3].$$
Its fixed point set is a complex geodesic and the $\mu$-reflection rotates
around this complex geodesic by the angle $\arg(\mu)$.
For $\mu=-1$ we get the usual complex reflections described above.
The complex $\mu$-reflection in a complex geodesic $C$ with a polar vector $c$
is given by
$$\iota_C^{\mu}(z)=z+(\mu-1)\cdot\frac{\<z,c\>}{\<c,c\>}\cdot c.$$

\myskip
{\bf Classification of isometries:}
An isometry $\varphi$ of $\chp$ is called {\it elliptic} if it has a fixed
point in $\chp$. It is called {\it hyperbolic\/} (or {\it loxodromic\/}) if
its {\it displacement} $d_{\varphi}=\inf\{d(x,\varphi(x))\st x\in\chp\}$ is positive.
Here $d$ is the inner metric on $\chp$.
A hyperbolic isometry has two fixed points that lie in $\dd\chp$.
An isometry, which is neither elliptic nor hyperbolic, is called {\it parabolic}.
It has one fixed point that lies in $\dd\chp$.
So far this was the usual classification of isometries of non-positive curved spaces.
In the case of the complex hyperbolic space we can refine this classification.
An elliptic element is called {\it regular elliptic\/} if all its eigenvalues
are distinct. A parabolic element, which can be written as an element of
$\U(2,1)$ with only eigenvalue~$1$, is called {\it unipotent}, otherwise the
parabolic element is called {\it ellipto-parabolic\/} (or {\it
skew-parabolic\/}).

\myskip
Using the discriminant function
$$\rho(z)=|z|^4-8\Re(z^3)+18|z|^2-27$$
we can classify isometries of the complex hyperbolic plane by the traces of the corresponding matrices:
An isometry $A\in\SU(2,1)$ is regular elliptic iff $\rho(\trace A)<0$
and hyperbolic iff $\rho(\trace A)>0$.
If $\rho(\trace A)=0$ there are three cases.
If $(\trace A)^3=27$ then $A$ is unipotent.
Otherwise, $A$ is either a complex reflection in a complex geodesic
or a complex reflection about a point,
or $A$ is ellipto-parabolic.
The proof can be found in~\cite{G99}, Theorem~6.2.4.
Note that for real $z$ the function~$\rho$ factors into $\rho(z)=(z+1)(z-3)^3$.
This means for $A\in\SU(2,1)$ whose trace is real, that $A$ is regular elliptic iff $\trace A\in(-1,3)$
and hyperbolic iff $\trace A\not\in[-1,3]$.

\section{Complex Hyperbolic Triangles and $\al$-Invariant} 

\label{cpxhyptri}

In this section we describe a parameterisation of the space of (complex
hyperbolic) triangles in~$\chp$,
\ie of triples $(C_1,C_2,C_3)$ of complex geodesics,
by means of an invariant $\al$.
Let $c_k$ be the normalised polar vector of the complex geodesic $C_k$.
Let $r_k=|\<c_{k-1},c_{k+1}\>|$.
If the complex geodesics $C_{k-1}$ and $C_{k+1}$ meet at the angle $\varphi_k$, then $r_k=\cos\varphi_k$.
If the complex geodesics $C_{k-1}$ and $C_{k+1}$ are ultra-parallel with distance $\ell_k$,
then $r_k=\cosh(\ell_k/2)$. We define the {\it angular invariant} $\al$ of the triangle $(C_1,C_2,C_3)$ by
$$\al=\arg\left(\prod_{k=1}^3\<c_{k-1},c_{k+1}\>\right).$$
The angular invariant is obviously invariant under isometries of $\chp$.
The complex reflection $\iota_k=\iota_{C_k}$ in the complex geodesic $C_k$ is defined by
$$\iota_k(z)=-z+2\<z,c_k\>\cdot c_k.$$
(Note that because of the property $\<c_k,c_k\>=1$ we can simplify the formula for the reflection.)

\begin{rem}
If the complex geodesics $C_{k-1}$ and $C_{k+1}$ or their closures intersect, 
then the vector $v_k=c_{k-1}\hcross c_{k+1}$ is a negative resp.\ null vector,
which corresponds to the intersection point,
\ie to a vertex of the complex hyperbolic triangle.
The points $[v_k]$ can be also meaningfully interpreted in the case of
ultra-parallel geodesics $C_{k-1}$ and $C_{k+1}$.
\end{rem}

We call a complex hyperbolic triangle $(C_1,C_2,C_3)$
a $(\varphi_1,\varphi_2,\varphi_3)$-triangle
if the complex geodesics $C_{k-1}$ and $C_{k+1}$ meet at the angle $\varphi_k$.

\begin{prop}
\label{ppp-triangle-existence}
A $(\varphi_1,\varphi_2,\varphi_3)$-triangle in $\chp$ is determined uniquely up to isometry
by the triple $(\varphi_1,\varphi_2,\varphi_3)$ and the angular invariant
$$\al=\arg\left(\prod_{k=1}^3\<c_{k-1},c_{k+1}\>\right).$$
For any $\al\in[0,2\pi]$ there exists a $(\varphi_1,\varphi_2,\varphi_3)$-triangle in $\chp$ with angular invariant~$\al$
if and only if
$$\cos\al<\frac{r_1^2+r_2^2+r_3^2-1}{2r_1r_2r_3},$$
where $r_k=\cos\varphi_k$.
\end{prop}

\begin{proof}
We can normalise so that
$$
  c_1=\begin{pmatrix}z\\\xi\\\be\end{pmatrix},\quad
  c_2=\begin{pmatrix}\ga\\\de\\0\end{pmatrix},\quad\text{and}\quad
  c_3=\begin{pmatrix}0\\1\\0\end{pmatrix}
$$
with $z\in\c$, $\xi,\be,\ga,\de\in\r$ and $\xi,\ga,\de>0$.
We compute
$$
  \<c_1,c_3\>=\xi,\quad
  \<c_2,c_1\>=\ga\bar z+\de\xi,\quad
  \<c_3,c_2\>=\de.
$$
The conditions $|\<c_{k-1},c_{k+1}\>|=r_k$ and $\<c_k,c_k\>=1$ imply
\begin{align*}
  \de&=r_1,\quad \ga=s_1,\\
  \xi&=r_2,\quad |z|>s_2,\quad \be=(|z|^2+r_2^2-1)^{1/2},\\
  z&=(r_3e^{-i\al}-r_1r_2)s_1^{-1},\\
\end{align*}
where $r_k=\cos(\varphi_k)$ and $s_k=\sin(\varphi_k)$.
The inequality $|z|>s_2$ finally implies $|r_3e^{i\al}-r_1r_2|>s_1s_2$.
Computation shows that this inequality is equivalent to
$\cos\al<(r_1^2+r_2^2+r_3^2-1)(2r_1r_2r_3)^{-1}$.
\end{proof}

\begin{rem}
In the case of angular invariant $\al=\pi$ all the vertices $[v_1]$, $[v_2]$, and $[v_3]$ of the triangle
lie in one real slice.
The corresponding triangle group representation stabilises this real slice.
The angular invariant $\al$ with
$$\cos\al=\frac{r_1^2+r_2^2+r_3^2-1}{2r_1r_2r_3}$$
corresponds to the case that all vertices $[v_1]$, $[v_2]$, and $[v_3]$ of the triangle coincide.
The groups with angular invariant $\al$ and $2\pi-\al$ are conjugate via an
anti-holomorphic isometry of~$\chp$.
For this reason we can restrict ourselves to the cases of angular invariant $\al\in(0,\pi]$.
\end{rem}

\myskip
Similar statements can be proved also in the cases when some of the complex
geodesics are ultra-parallel.
For example, for an ultra-parallel \lll-triangle $(C_1,C_2,C_3)$, \ie
in the case that the complex geodesics $C_{k-1}$ and $C_{k+1}$ are ultra-parallel with distance~$\ell_k$,
we obtain by similar reasoning the following result:

\begin{prop}
\label{lll-triangle-existence}
An ultra-parallel \lll-triangle in $\chp$ is determined
uniquely up to isometry by the triple $(\ell_1,\ell_2,\ell_3)$ and the angular invariant
$$\al=\arg\left(\prod_{k=1}^3\<c_{k-1},c_{k+1}\>\right).$$
For any $\al\in[0,2\pi]$ there exists an ultra-parallel \lll-triangle in $\chp$ with angular invariant $\al$
if and only if
$$\cos\al<\frac{r_1^2+r_2^2+r_3^2-1}{2r_1r_2r_3},$$
where $r_k=\cosh(\ell_k/2)$.
\end{prop}

{\bf Examples of explicit parameterisations:}
In the case $p_3=\infty$ we consider for $\al\in[-\pi,\pi]$
the complex geodesics with normalised polar vectors
$$
  c_1=\begin{pmatrix}1\\z_2\\-z_2\end{pmatrix},\quad
  c_2=\begin{pmatrix}1\\z_1\\-z_1\end{pmatrix},\quad\text{and}\quad
  c_3=\begin{pmatrix}0\\1\\0\end{pmatrix},
$$
where $z_1=\cos(\pi/p_1)e^{-i\al/2}$ and $z_2=\cos(\pi/p_2)e^{i\al/2}$.
Computing
$$
  \<c_1,c_3\>=r_2e^{i\al/2},\quad
  \<c_2,c_1\>=1,\quad
  \<c_3,c_2\>=r_1e^{i\al/2}
$$
we see that $\ga_k\mapsto\iota_{C_k}$ defines a $(p_1,p_2,\infty)$-representation with $\al$-invariant equal to $\al$.
For the sake of completeness we also compute the vertices of the triangle
$$
  [v_1]=[-r_1e^{i\al/2}:0:1],\quad
  [v_2]=[-r_2e^{-i\al/2}:0:1],\quad
  [v_3]=[0:1:-1].
$$
A special case of this parameterisation is the parameterisation of \ppi-groups in~\cite{W}.
Another example of explicit parameterisation is the parameterisation of \uuu-groups in~\cite{Sch2003}.

\section{Connection of the Angular Invariant $\al$ and Other Invariants}

\label{anglinvcompar}

In this section we compare our parameterisation of the space of complex
hyperbolic triangles with other parameterisations by Ulrich Brehm~\cite{B}
and Jeffrey Hakim and Hanna Sandler~\cite{HS:Bruhat}.

Let $(C_1,C_2,C_3)$ be a complex hyperbolic triangle.
Let $c_k$ be the normalised polar vector of the complex geodesic $C_k$.
Assume that for $k=1,2,3$ the complex geodesics $C_{k-1}$ and $C_{k+1}$ or their closures intersect.
Let $v_k=c_{k-1}\hcross c_{k+1}$.

\myskip
The {\it Cartan angular invariant}
of three points $ [v_1]$, $[v_2]$, and $[v_3]$ in $\dd\chp$
was defined by \'E.~Cartan in~\cite{C32} as follows:
$$A([v_1],[v_2],[v_3])=\arg(-\<v_1,v_2\>\<v_2,v_3\>\<v_3,v_1\>).$$
To compare the Cartan angular invariant and our invariant $\al$ we need the following proposition:

\begin{prop}
\label{hcross-formulas}
Let $\<c_{k-1},c_{k+1}\>=r_ke^{i\theta_k}$.
Then we have $\<v_k,v_k\>=r_k^2-1$,
$\<v_{k-1},v_{k+1}\>=e^{i\theta_k}(r_k-r_{k-1}r_{k+1}e^{-i\al})$, and
$$
  \prod\limits_{k=1}^3\<v_k,v_{k+1}\>
  =\prod\limits_{k=1}^3\overline{\<v_{k-1},v_{k+1}\>}
  =e^{-i\al}\cdot\prod\limits_{k=1}^3(r_k-r_{k-1}r_{k+1}e^{i\al}).
$$
\end{prop}

\begin{proof}
The proof of this proposition, based on the fact that $[v_k]=[c_{k-1}\hcross c_{k+1}]$
and the formula for $\<a\hcross c,b\hcross c\>$, is straightforward.
\end{proof}

Using the last proposition we obtain the following formula for the Cartan angular invariant
$$A=\arg(-e^{-i\al}(1-e^{i\al})^3)=\frac{\al-\pi}{2}\moddpi$$
because of $\arg(1-e^{i\al})=(\al-\pi)/2$.

Ulrich Brehm defined in~\cite{B} the {\it shape invariant\/}
of three points $ [v_1]$, $[v_2]$, and $[v_3]$ in $\chp$ as
$$\si([v_1],[v_2],[v_3])=\Re\left(\frac{\<v_1,v_2\>\<v_2,v_3\>\<v_3,v_1\>}{\<v_1,v_1\>\<v_2,v_2\>\<v_3,v_3\>}\right).$$
The shape invariant $\si$ is invariant under $\PU(2,1)$,
and the triangle $[v_1]$, $[v_2]$, $[v_3]$ is determined uniquely up to
isometry by its three angles together with the shape invariant.
Using proposition~\ref{hcross-formulas} we obtain the following formula that describes
the relation between the Brehm's shape invariant $\si$ and our invariant $\al$:
$$
  \si
  =\frac{r_1^2r_2^2r_3^2\cos2\al-r_1r_2r_3(r_1^2+r_2^2+r_3^2+1)\cos\al+(r_1^2r_2^2+r_2^2r_3^2+r_3^2r_1^2)}
        {(1-r_1^2)(1-r_2^2)(1-r_3^2)}.
$$
In \cite{HS:Bruhat} Jeffrey Hakim and Hanna Sandler defined such angular
invariants also in the cases of triangles with one or two ideal vertices:
$$\eta([v_1],[v_2],[v_3])=\frac{\<v_3,v_1\>\<v_1,v_2\>}{\<v_3,v_2\>\<v_1,v_1\>}.$$
Using proposition~\ref{hcross-formulas} we compute
$$\eta=e^{-i\al}\cdot\frac{(r_2-r_1e^{i\al})(1-r_1r_2e^{i\al})}{(r_1-r_2e^{i\al})(r_1^2-1)}.$$

\section{Winding Numbers and Other Combinatorial Functions}

\label{winding}

In this section we introduce some combinatorial functions of words,
which will appear in the combinatorial trace formula.

\myskip
We first fix some notation. Elements of the \ppp-triangle group, or the abstract group $\Gappp$,
may be naturally identified with finite sequences of numbers $1$, $2$, and $3$.
We refer to such sequences as {\defit words}.
The sequence $a=(a_1,\dots,a_n)$ corresponds to the element
$\iota_a:=\iota_{a_1}\cdots\iota_{a_n}$.
We shall simplify notation where it is convenient,
for example we write $(123)$ instead of $(1,2,3)$
and $\iota_{123}$ instead $\iota_{(123)}$ or $\iota_{(1,2,3)}$.
The empty sequence corresponds to the identity element denoted by $\iota_{\emptyset}$.
There is a natural definition of the powers $(a_1,\dots,a_n)^w$ for $w\in\z$,
which is consistent with taking powers of the corresponding elements in the group,
for example $(123)^2=(123123)$ and $(123)^{-1}=(321)$.
We denote by $\tau_a$ the trace of the element $\iota_a$.

\myskip
We now introduce the winding number.
Let $\chi$ denote the unique nontrivial character modulo~$3$, \ie if
$a$ is any integer then $\chi(a)$ is the unique element
of the set of integers $\{-1,0,1\}$,
which is congruent to~$a$ modulo~$3$.
Define
$$w(a_1,\dots,a_n)=\frac13\cdot\sum\limits_{m=1}^n\,\chi(a_{m+1}-a_m),$$
where $a_{n+1}:=a_1$.

\begin{rem}
The integer $w(a_1,\dots,a_n)$ can be interpreted as the winding number of the
loop $a_1\to a_2\to\cdots\to a_n\to a_1$ along the vertices of the $
(1,2,3)$-triangle (see figure~\ref{figa}).


\begin{figure}
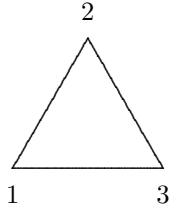

  \begin{center}
    \forcehmode
      \bgroup
        \beginpicture
          \setcoordinatesystem units <2cm,2cm>
          \plot 1 0 1.5 0.866 2 0 1 0 /
          \put {$1$} [] <0pt,-10pt> at 1 0
          \put {$2$} [] <0pt,10pt> at 1.5 0.866
          \put {$3$} [] <0pt,-10pt> at 2 0
        \endpicture
      \egroup
  \end{center}
  \caption{Winding Number and $(1,2,3)$-Triangle}
  \label{figa}
\end{figure}
\end{rem}

We are interested in traces of elements.
For $A_1,\dots,A_m\in\SU(2,1)$ we have
$$\trace(A_1\cdots A_m)=\trace(A_m\cdot A_1\cdots A_{m-1}).$$
Since the trace of a product does not change under cyclic permutation of the factors,
we consider cyclic words instead of usual (linear) words.
A {\defit cyclic word\/} is the orbit of a linear word under cyclic permutations.

\myskip
Dealing with winding numbers we often use the following consideration:
We start with the following operations on cyclic words:
reduction
$$(\cdots,k,k,\cdots)\rightarrow(\cdots,k,\cdots)$$
and straightening
$$(\cdots,k,l,k,\cdots)\rightarrow(\cdots,k,\cdots).$$
These operations do not change the winding number,
and any cyclic word $a$ can be transformed into the cyclic word $(123)^{w(a)}$
using these operations.

\myskip
We now introduce combinatorial functions $u_k$
that will be used in the combinatorial trace formula.

\myskip
Let
$$u_k(a_1,\dots,a_n)=|\{m\in\{1,\dots,n\}\st\{a_m,a_{m+1}\}=\{k-1,k+1\}\}|,$$
where $a_{n+1}=a_1$.
Using the function $\psi_k$ given by $\psi_k(a,b)=1$ iff
$$\{a,b\}=\{k-1,k+1\}$$
and $\psi_k(a,b)=0$ otherwise we can describe $u_k$ as
$$u_k(a_1,\dots,a_n)=\sum\limits_{m=1}^n\,\psi_k(a_m,a_{m+1}),$$
where again $a_{n+1}=a_1$.

\section{Combinatorial Trace Formula}

\label{combtraceformula}

The first main result of this paper is the following combinatorial formula for
traces of elements in complex hyperbolic triangle group.

\myskip
We first fix some notation.
For a word $a=(a_1,\dots,a_n)$ and a subset $S$ of the set $[n]:=\{1,\dots,n\}$
we denote by $u_k(S)$ resp.\ $w(S)$ the values of $u_k$ resp.\ $w$
for the corresponding subsequences of~$a$, for example
$u_k(S)=u_k(a_{i_1},\dots,a_{i_m})$ for $S=\{i_1,\dots,i_m\}$ with $1\le i_1<\dots<i_m\le n$.

\begin{thm}
\label{comb}
Let $a=(a_1,\dots,a_n)$ be a cyclic word, $\iota_a=\iota_{a_1}\cdots\iota_{a_n}$ the
corresponding element of the \ppp-group and $\tau_a$ the trace of $\iota_a$.
Then
$$
  \tau_a
  =(-1)^n\cdot
   \left(
         2+
         \sum\limits_{S\subset[n]}
         (-2)^{|S|}\cdot r_1^{u_1(S)}r_2^{u_2(S)}r_3^{u_3(S)}\cdot e^{i\al w(S)}
   \right).
$$
Equivalently, $\tau_a=(-1)^n\cdot(2+p_a(e^{i\al}))$,
where $p_a(z)=\sum_w q_w z^w$ is the Fourier polynomial with coefficients
$$
  q_w=\sum\limits_{S\subset[n]\atop w(S)=w}
      (-2)^{|S|}\cdot r_1^{u_1(S)}r_2^{u_2(S)}r_3^{u_3(S)}.
$$
\end{thm}

\begin{proof}
We write the complex reflections as $\iota_k=-\id+2c_kc_k^*$,
where $c^*$ is the dual vector to a vector $c\in\cdo$, \ie $c^*(z)=\<z,c\>$.
To get a formula for traces we expand
\begin{align*}
  \iota_a
  &=\iota_{a_1}\cdots\iota_{a_n}\\
  &=(-\id+2c_{a_1}c_{a_1}^*)\cdots(-\id+2c_{a_n}c_{a_n}^*)\\
  &=\sum\limits_{S=\{i_1,\dots,i_m\}\subset[n]}(-1)^{n-|S|}\cdot2^{|S|}\cdot(c_{a_{i_1}}c_{a_{i_1}}^*\cdots c_{a_{i_m}}c_{a_{i_m}}^*).\\
\end{align*}
Let $S=\{i_1,\dots,i_m\}$ be a subset of $[n]$.
For ease of notation let $b_j=a_{i_j}$ and $b=(b_1,\dots,b_m)$.
We have
$$
  c_{b_1}c_{b_1}^*c_{b_2}c_{b_2}^*\cdots c_{b_{m-1}}c_{b_{m-1}}^*c_{b_m}c_{b_m}^*
  =c_{b_1}^*(c_{b_2})\cdots c_{b_{m-1}}^*(c_{b_m})\cdot(c_{b_1}c_{b_m}^*)
$$
and $\trace(c_{b_1}c_{b_m}^*)=c_{b_m}^*(c_{b_1})$ implies
\begin{align*}
  \trace(c_{b_1}c_{b_1}^*\cdots c_{b_m}c_{b_m}^*)
  &=c_{b_1}^*(c_{b_2})\cdots c_{b_{m-1}}^*(c_{b_m})\cdot c_{b_m}^*(c_{b_1})\\
  &=\<c_{b_2},c_{b_1}\>\cdots\<c_{b_m},c_{b_{m-1}}\>\cdot\<c_{b_1},c_{b_m}\>.
\end{align*}
We claim that
$$
  \<c_{b_2},c_{b_1}\>\cdots\<c_{b_m},c_{b_{m-1}}\>\cdot\<c_{b_1},c_{b_m}\>
  =r_1^{u_1(b)}r_2^{u_2(b)}r_3^{u_3(b)}\cdot e^{i\al w(b)}.
$$
It is clear that the absolute value of the complex number on the left hand
side of the equation is equal to $r_1^{u_1(b)}r_2^{u_2(b)}r_3^{u_3(b)}$.
The statement about the argument of this complex number is clear for
the words of the form $(123)^w$ for some $w\in\z$.
Any other cyclic word $b$ can be transformed into the word $(123)^{w(b)}$
using the winding number preserving operations
of reduction and straightening as described in section~\ref{winding}.
It remains to note that this operations do not change the argument of the
complex number on the left hand side of the equation.
So we have proven the formula
$$
  \<c_{a_{i_2}},c_{a_{i_1}}\>\cdots\<c_{a_{i_m}},c_{a_{i_{m-1}}}\>\cdot\<c_{a_{i_1}},c_{a_{i_m}}\>
  =r_1^{u_1(S)}r_2^{u_2(S)}r_3^{u_3(S)}\cdot e^{i\al w(S)}
$$
for $S=\{i_1,\dots,i_m\}\subset[n]$.
(Note that this formula does not hold for $S=\emptyset$, so we have to handle
this case separately.)
Using this formula we finally obtain
$$
  \tau_a
  =(-1)^n\cdot
   \left(
         2+
       \sum\limits_{S\subset[n]}
       (-2)^{|S|}\cdot r_1^{u_1(S)}r_2^{u_2(S)}r_3^{u_3(S)}\cdot e^{i\al w(S)}
   \right).\qedhere
$$
\end{proof}

\begin{rem}
In the ideal triangle group case, \ie for $p_1=p_2=p_3=\infty$, it follows from the
Binomial Theorem
$$\sum\limits_{k}q_k=\sum\limits_{S\subset[n]}(-2)^{|S|}=(1+(-2))^n=(-1)^n.$$
\end{rem}

\section{Recursive Trace Formula}

\label{rectraceformula}

In explicit computations we rather use the following recursive trace formula,
which can be derived from the combinatorial trace formula of the last section.
We first state some properties of the winding number and of the combinatorial function $u_k$.
They follow from routine induction arguments.

\begin{prop}
\label{winding-formulas}
We note that
$$w(a_1,\dots,a_n)=w(a_1,\dots,a_m)+w(a_m,\dots,a_n)+w(a_1,a_m,a_{m+1},a_n)$$
for $m=2,\dots,n-1$.
Especially
$$w(a_1,\dots,a_n,a_{n+1})=w(a_1,\dots,a_n)+w(a_1,a_n,a_{n+1}),$$
and
$$w(b_1,\dots,b_m,a_{n-2},a_{n-1},a_n)=w(b_1,\dots,b_m,a_{n-2},a_n)+w(a_{n-2},a_{n-1},a_n).\qedhere$$
\end{prop}

\begin{prop}
\label{uk-formulas}
We note that
$$u_k(a_1,\dots,a_n,a_{n+1})=u_k(a_1,\dots,a_n)+v_k(a_n,a_{n+1},a_1)$$
and
$$u_k(b_1,\dots,b_k,a_{n-2},a_{n-1},a_n)=u_k(b_1,\dots,b_k,a_{n-2},a_n)+v_k(a_{n-2},a_{n-1},a_n),$$
where
$$v_k(a,b,c)=\psi_k(a,b)+\psi_k(b,c)-\psi_k(a,c).$$
\end{prop}

Define {\defit deletion operators\/} $\de$, $\de'$, and $\de''$:
For a word $a=(a_1,\dots,a_n)$ let
\begin{align*}
  \de&(a)=(a_1,\dots,a_{n-3},a_{n-2},a_{n-1}),\\
  \de'&(a)=(a_1,\dots,a_{n-3},a_{n-2},a_n),\\
  \de''&(a)=(a_1,\dots,a_{n-3},a_{n-1},a_n).
\end{align*}

\begin{thm}
\label{reccur}
Let $a=(a_1,\dots,a_n)$ be a word of length $n\ge3$
and $\tau_a$ the trace of the corresponding element $\iota_a$.
Then
$$
  \tau_a=-(\tau_{\de a}+\tau_{\de''a}+\tau_{\de\de''a})
         +\be(\tau_{\de'a}+\tau_{\de\de'a}+\tau_{\de'\de''a}+\tau_{\de\de'\de''a}),
$$
where
$$\be=2r_1^{v_1(T)}r_2^{v_2(T)}r_3^{v_3(T)}\cdot e^{i\al w(T)}-1$$
with $T=(a_{n-2},a_{n-1},a_n)$.
\end{thm}

\begin{proof}
The proof is along the lines of the proof of Theorem~3 in~\cite{S95}.
We consider
\begin{align*}
  t_a:&=(-1)^n\tau_a-2\\
      &=\sum\limits_{S\subset[n]}
        (-2)^{|S|}\cdot r_1^{u_1(S)}r_2^{u_2(S)}r_3^{u_3(S)}e^{i\al w(S)},\\
  t'_a:&=t_a-t_{\de a}
        =(-1)^n(\tau_a+\tau_{\de a})\\
       &=\sum\limits_{n\in S\subset[n]}
         (-2)^{|S|}\cdot r_1^{u_1(S)}r_2^{u_2(S)}r_3^{u_3(S)}e^{i\al w(S)},\\
  t''_a:&=t'_a-t'_{\de' a}
         =(-1)^n(\tau_a+\tau_{\de a}+\tau_{\de'a}+\tau_{\de\de'a})\\
        &=\sum\limits_{n-1,n\in S\subset[n]}
          (-2)^{|S|}\cdot r_1^{u_1(S)}r_2^{u_2(S)}r_3^{u_3(S)}e^{i\al w(S)},\\
  t'''_a:&=t''_a-t''_{\de'' a}\\
         &=(-1)^n(\tau_a+\tau_{\de a}+\tau_{\de'a}+\tau_{\de\de'a}+\tau_{\de''a}+\tau_{\de\de''a}+\tau_{\de'\de''a}+\tau_{\de\de'\de''a})\\
         &=\sum\limits_{n-2,n-1,n\in S\subset[n]}
           (-2)^{|S|}\cdot r_1^{u_1(S)}r_2^{u_2(S)}r_3^{u_3(S)}e^{i\al w(S)}.\\
\end{align*}
For $S\subset[n-3]$ let $\hat S:=S\cup\{n-2,n-1,n\}$ and $S':=S\cup\{n-2,n\}$.
By propositions~\ref{winding-formulas} and \ref{uk-formulas} it holds for $S\subset[n-3]$
$$
  w(\hat S)=w(S')+w(T)
  \quad\hbox{and}\quad
  u_k(\hat S)=u_k(S')+v_k(T).
$$
This implies
\begin{align*}
  t'''_a
  &=\sum\limits_{n-2,n-1,n\in S\subset[n]}
    (-2)^{|S|}\cdot r_1^{u_1(S)}r_2^{u_2(S)}r_3^{u_3(S)}e^{i\al w(S)}\\
  &=\sum\limits_{S\subset[n-3]}
    (-2)^{|S|+3}\cdot r_1^{u_1(\hat S)}r_2^{u_2(\hat S)}r_3^{u_3(\hat S)}e^{i\al w(\hat S)}\\
  &=-2\cdot r_1^{v_1(T)}r_2^{v_2(T)}r_3^{v_3(T)}e^{i\al w(T)}\cdot
    \sum\limits_{S\subset[n-3]}
    (-2)^{|S|+2}\cdot r_1^{u_1(S')}r_2^{u_2(S')}r_3^{u_3(S')}e^{i\al w(S')}\\
  &=-(\be+1)\cdot t''_{\de'a}.\\
\end{align*}
From this equation and the definition of $t'''_a$ we obtain
$$
  t''_a
  =t''_{\de'' a}+t'''_a
  =t''_{\de'' a}-(\be+1)\cdot t''_{\de'a}.
$$
On the other hand $t''$ can be expressed in terms of the traces as
\begin{align*}
  &(-1)^n(\tau_a+\tau_{\de a}+\tau_{\de'a}+\tau_{\de\de'a})\\
  &\qquad=(-1)^{n-1}(\tau_{\de''a}+\tau_{\de\de''a}+\tau_{\de'\de''a}+\tau_{\de\de'\de''a})\\
  &\qquad\phantom{=}-(\be+1)\cdot(-1)^{n-1}(\tau_{\de'a}+\tau_{\de\de'a}+\tau_{\de'\de'a}+\tau_{\de\de'\de'a}).\\
\end{align*}
This implies using $\de'\de'=\de'\de''$
\begin{align*}
  \tau_a
  &=-(\tau_{\de''a}+\tau_{\de\de''a}+\tau_{\de'\de''a}+\tau_{\de\de'\de''a}+\tau_{\de a}+\tau_{\de'a}+\tau_{\de\de'a})\\
  &\phantom{=}+(\be+1)\cdot(\tau_{\de'a}+\tau_{\de\de'a}+\tau_{\de'\de''a}+\tau_{\de\de'\de''a})\\
  &=-(\tau_{\de''a}+\tau_{\de\de''a}+\tau_{\de a})+\be\cdot(\tau_{\de'a}+\tau_{\de\de'a}+\tau_{\de'\de''a}+\tau_{\de\de'\de''a}).\qedhere\\
\end{align*}
\end{proof}

In order to use this formula to compute traces recursively, it is necessary to
know the traces of all elements associated to words of length $0$, $1$, and $2$:
\begin{align*}
  &\tau_{\emptyset}=3,\quad
  \tau_k=-1,\\
  &\tau_{k-1,k+1}=\tau_{k+1,k-1}=4r_k^2-1.\\
\end{align*}

\section{Examples}

\label{examples}

We compute the traces for some special elements of complex hyperbolic triangle
groups, in particular for $\wA=\iota_{1323}$ and $\wB=\iota_{123}$.

\myskip\noindent
We have
$$\tau_{123}=-(\tau_{12}+\tau_{23}+\tau_2)+\beta\cdot(\tau_{13}+\tau_1+\tau_3+\tau_{\emptyset}),$$
where
$$\beta=2r_1^{v_1(1,2,3)}r_2^{v_2(1,2,3)}r_3^{v_3(1,2,3)}e^{i\al w(1,2,3)}-1=2r_1r_2^{-1}r_3e^{i\al}-1.$$
Hence
\begin{align*}
  \tau_{123}&=-((4r_3^2-1)+(4r_1^2-1)+(-1))+(2r_1r_2^{-1}r_3e^{i\al}-1)\cdot((4r_2^2-1)+1)\\
            &=8r_1r_2r_3e^{i\al}-(4(r_1^2+r_2^2+r_3^2)-3).\\
\end{align*}

\myskip\noindent
Furthermore we have
$$\tau_{2321}=-(\tau_{232}+\tau_{221}+\tau_{22})+\beta\cdot(\tau_{231}+\tau_{23}+\tau_{21}+\tau_2),$$
where
$$\beta=2r_1^{v_1(3,2,1)}r_2^{v_2(3,2,1)}r_3^{v_3(3,2,1)}e^{i\al w(3,2,1)}-1=2r_1r_2^{-1}r_3e^{-i\al}-1.$$
It holds $\iota_{231}=\iota_1(\iota_{123})\iota_1^{-1}$ and therefore $\tau_{231}=\tau_{123}$.
Because of $\tau_{232}+\tau_{221}+\tau_{22}=\tau_3+\tau_1+\tau_{\emptyset}=1$ and
\begin{align*}
  &\tau_{231}+\tau_{23}+\tau_{21}+\tau_2\\
  &=(8r_1r_2r_3e^{i\al}-(4r_1^2+4r_2^2+4r_3^2-3))+(4r_1^2-1)+(4r_3^2-1)+(-1)\\
  &=8r_1r_2r_3e^{i\al}-4r_2^2=4r_2(2r_1r_3e^{i\al}-r_2)\\
\end{align*}
we obtain
\begin{align*}
  \tau_{2321}&=-1+(2r_1r_2^{-1}r_3e^{-i\al}-1)\cdot4r_2(2r_1r_3e^{i\al}-r_2)\\
             &=4(2r_1r_3e^{-i\al}-r_2)(2r_1r_3e^{i\al}-r_2)-1\\
             &=4\cdot|2r_1r_3e^{i\al}-r_2|^2-1\\
             &=(16r_1^2r_3^2+4r_2^2-1)-16r_1r_2r_3\cos\al.\\
\end{align*}

\myskip\noindent
Similarly, for any $k$ we get
\begin{align*}
  \si_k
  &:=\tau_{k,k-1,k,k+1}
  =\tau_{k,k+1,k,k-1}\\
  &=4\cdot|2r_{k-1}r_{k+1}e^{i\al}-r_k|^2-1\\
  &=(16r_{k-1}^2r_{k+1}^2+4r_k^2-1)-16r_1r_2r_3\cos\al.\\
\end{align*}

\section{Some Properties of the Traces}

\label{someproperties}

In this section we study some properties of our trace formulas.
Lemma~\ref{qw-in-ring} will be also used later in section~\ref{appl-disc}.

\begin{lem}
\label{qw-in-ring}
The coefficients $q_w$ in the combinatorial trace formula satisfy
$$q_w\in(8r_1r_2r_3)^{|w|}\cdot\z[4r_1^2,4r_2^2,4r_3^2].$$
\end{lem}

\begin{proof}
The coefficient $q_w$ is a sum of the terms of the form
$$\pm 2^{|b|}\cdot r_1^{u_1(b)}\cdot r_2^{u_2(b)}\cdot r_3^{u_3(b)}$$
for some word $b$ of length~$|b|$ and winding number~$w(b)=w$.
For the words $b=(123)^w$ with $w\in\z$
we have
$$
  2^{|b|}\cdot r_1^{u_1(b)}\cdot r_2^{u_2(b)}\cdot r_3^{u_3(b)}
  =2^{3|w|}\cdot r_1^{|w|}\cdot r_2^{|w|}\cdot r_3^{|w|}
  =(8r_1r_2r_3)^{|w|}.
$$
Any other word $b$ can be transformed into the word $(123)^{w(b)}$
using the winding number preserving operations described in section~\ref{winding}.
But the number
$$\pm 2^{|b|}\cdot r_1^{u_1(b)}\cdot r_2^{u_2(b)}\cdot r_3^{u_3(b)}$$
is divided by $2$ under reduction
and by $4r_k^2$ for some $k\in\{1,2,3\}$ under straightening.
Altogether this implies
$$
  \pm 2^{|b|}\cdot r_1^{u_1(b)}\cdot r_2^{u_2(b)}\cdot r_3^{u_3(b)}
  \in(8r_1r_2r_3)^{|w(b)|}\cdot\z[4r_1^2,4r_2^2,4r_3^2]
$$
and hence
$q_w\in(8r_1r_2r_3)^{|w|}\cdot\z[4r_1^2,4r_2^2,4r_3^2]$.
\end{proof}

\begin{thm}
For any cyclic word~$a$ of length~$n$ there is a polynomial $s_a(z)$ in
$\z[4r_1^2,4r_2^2,4r_3^2][z,\bz]$ of degree at most $n/2$ in $4r_k^2$ and at most $n/3$ in $z$ and $\bz$,
independent of the angular invariant $\al$, such that
$\tau_a=s_a(\tau)$, where $\tau=\tau_{123}$.
\end{thm}

\begin{rem}
In particular $s_a=B_k$ for $a=(123)^k$, where polynomials $B_k$ are the
analogues of the Chebyshev polynomials for $\SU(2,1)$ as defined by Hanna Sandler in~\cite{S95}.
\end{rem}

\begin{proof}
Let $R:=\z[4r_1^2,4r_2^2,4r_3^2]$.
The trace $\tau=\tau_{123}$ is equal to $\tau=(8r_1r_2r_3)\cdot e^{i\al}-c$,
where $c=4r_1^2+4r_2^2+4r_3^2-3\in R$.
Therefore
$$
  e^{i\al w}
  =(8r_1r_2r_3)^{-|w|}\cdot([\tau~\hbox{or}~\bar\tau]+c)^w
  \in(8r_1r_2r_3)^{-|w|}\cdot R[\tau,\bar\tau].
$$
On the other hand the coefficients $q_w$ in the combinatorial trace formula
satisfy $q_w\in(8r_1r_2r_3)^{|w|}\cdot R$ by lemma~\ref{qw-in-ring},
hence $q_w\cdot e^{i\al w}\in R[\tau,\bar\tau]$.
\end{proof}

\section{Groups Generated by Complex $\mu$-Reflections}

\label{mu-refl}

In this section we generalise the trace formulas for the case of groups
generated by complex $\mu$-reflections.

\myskip
Let $(C_1,C_2,C_3)$ be a complex hyperbolic triangle.
Let $c_k$ be the normalised polar vector of the complex geodesic $C_k$,
and let $r_k=|\<c_{k-1},c_{k+1}\>|$.
The angular invariant $\al$ of the triangle $(C_1,C_2,C_3)$ is
$$\al=\arg\left(\prod_{k=1}^3\<c_{k-1},c_{k+1}\>\right).$$
Let $\mu_1$, $\mu_2$, and $\mu_3$ be unit complex numbers.
The complex $\mu_k$-reflection $\iota_k=\iota_{C_k}^{\mu_k}$ in the complex geodesic $C_k$ is given by
$$\iota_k(z)=z+(\mu_k-1)\cdot\<z,c_k\>\cdot c_k.$$
We consider the subgroup generated by the reflections $\iota_1$, $\iota_2$, and $\iota_3$.
The formulas for traces can be generalised for this case.

\myskip
Let
$$n_k(a_1,\dots,a_n)=|\{m\in\{1,\dots,n\}\st a_m=k\}|.$$

\myskip
For a word $a=(a_1,\dots,a_n)$ and a subset $S$ of the set $[n]:=\{1,\dots,n\}$
we denote by $n_k(S)$, $u_k(S)$, resp.~$w(S)$ the values of $n_k$, $u_k$, resp.~$w$
for the corresponding subsequences of~$a$, for example
$u_k(S)=u_k(a_{i_1},\dots,a_{i_m})$ for $S=\{i_1,\dots,i_m\}$ with $1\le i_1<\dots<i_m\le n$.

\begin{thm}
\label{mu-comb}
Let $a=(a_1,\dots,a_n)$ be a word, $\iota_a=\iota_{a_1}\cdots\iota_{a_n}$ the
corresponding element of the group generated by the three reflections $\iota_k=\iota_{C_k}^{\mu_k}$
and $\tau_a$ the trace of~$\iota_a$.
Then
$$
  \tau_a
  =2
  +\sum\limits_{S\subset[n]}
     \left(\prod\limits_{k=1}^3(\mu_k-1)^{n_k(S)}\right)\cdot
     \left(\prod\limits_{k=1}^3 r_k^{u_k(S)}\right)\cdot
     e^{i\al w(S)}.
$$
\end{thm}

\begin{proof}
We write the complex $\mu_k$-reflection as $\iota_k=\id+(\mu_k-1)c_kc_k^*$,
where $c^*$ is the dual vector to a vector $c\in\cdo$, \ie $c^*(z)=\<z,c\>$.
In order to get a formula for traces we expand
\begin{align*}
  \iota_a
  &=\iota_{a_1}\cdots\iota_{a_n}\\
  &=(\id+(\mu_{a_1}-1)c_{a_1}c_{a_1}^*)\cdots(\id+(\mu_{a_n}-1)c_{a_n}c_{a_n}^*)\\
  &=\sum\limits_{S=\{i_1,\dots,i_m\}\subset[n]}
        \left(\prod\limits_{k=1}^3(\mu_k-1)^{n_k(S)}\right)\cdot
        (c_{a_{i_1}}c_{a_{i_1}}^*\cdots c_{a_{i_m}}c_{a_{i_m}}^*).\\
\end{align*}
Now we proceed as in the proof of theorem~\ref{comb}.
\end{proof}

The generalisations of other statements based on the combinatorial trace formula
for the case of the groups generated by $\mu$-reflections are straightforward,
for instance the generalisation of the recursive trace formula and of lemma~\ref{qw-in-ring}.

\myskip
As an example, we consider Mostow's non-arithmetic complex hyperbolic lattices~\cite{M}.
Let us denote by $\Gamma(p,\rho)$ a complex hyperbolic triangle groups generated by complex
$\mu$-reflections of the same finite order $p\in\{3,4,5\}$
$$\mu:=\mu_1=\mu_2=\mu_3=e^{\frac{2\pi i}{p}}$$
in the sides of an equiangular triangle with
$$r:=r_1=r_2=r_3=\frac{1}{2\sin\frac{\pi}{p}}$$
and the value of the angular invariant $\al$ of the form
$$\al=\frac{2\pi}{\rho}+\frac{\pi}{p}-\frac{\pi}{2}.$$

\begin{prop}
The traces of elements in the group $\Gamma(p,\rho)$ satisfy
$$\tau_a\in\q[e^{\frac{2\pi i}{p}},e^{\frac{2\pi i}{\rho}}].$$
\end{prop}

\begin{proof}
Generalising lemma~\ref{qw-in-ring} we obtain that the coefficients
$$
  q_w=\sum\limits_{S\subset[n]\atop w(S)=w}
      \left(\prod\limits_{k=1}^3(\mu_k-1)^{n_k(S)}\right)\cdot
      \left(\prod\limits_{k=1}^3 r_k^{u_k(S)}\right)
$$
in the combinatorial trace formula satisfy
$$q_w\in((\mu-1)r)^{3|w|}\cdot\z[\mu,(\mu-1)^2r^2].$$
For $\mu=e^{\frac{2\pi i}{p}}$ and $r=\frac{1}{2\sin\frac{\pi}{p}}$ it holds
$$(\mu-1)r=i\cdot e^{\frac{\pi i}{p}},$$
hence we obtain 
$$q_w\in(i e^{\frac{\pi i}{p}})^{3|w|}\cdot\z[e^{\frac{2\pi i}{p}}]$$
and therefore 
$$\tau_a\in\q[e^{\frac{2\pi i}{p}},ie^{\frac{\pi i}{p}}e^{\al}].$$
Finally 
$$i\cdot e^{-\frac{\pi i}{p}}\cdot e^{i\al}=e^{\frac{2\pi i}{\rho}}$$
implies, that the traces of elements in the group $\Gamma(p,\rho)$ satisfy
$$\tau_a\in\q[e^{\frac{2\pi i}{p}},e^{\frac{2\pi i}{\rho}}].\qedhere$$
\end{proof}

Mostow proves this result in lemma~17.2.1 in~\cite{M} using 
implicitely a trace formula similar to our combinatorial trace formula
for the special case of the groups $\Gamma(p,\rho)$.

\section{Applications: Non-Discreteness of Triangle Groups}

\label{appl-non-disc}

In this section we describe some necessary conditions for a triangle group representation
to be a discrete embedding.
For the computations in this and the next section the parameterisation of
complex hyperbolic \ppp-triangle groups by
$$t=\ctan\frac{\al}{2}=\sqrt{\frac{1+\cos\al}{1-\cos\al}}$$
is more suitable than the parameterisation by the angular invariant $\al$ itself. 
Recall that a \ppp-triangle group exists for $t\in(-\tinf,\tinf)$,
where
$$
  \tinf=\sqrt{\frac{1+\cinf}{1-\cinf}}
  \quad\hbox{and}\quad
  \cinf=\frac{r_1^2+r_2^2+r_3^2-1}{2r_1r_2r_3}.
$$
(For $\cinf>1$ resp.\ $\cinf<-1$ we set $\tinf=+\infty$ resp.\ $\tinf=-\infty$.)

\myskip
Let $\phi:\Gappp\to\PU(2,1)$ be a complex hyperbolic triangle group representation
and $G:=\phi(\Gappp)$ the corresponding complex hyperbolic triangle group.
Assume that $\ga$ is an element of infinite order in $\Gappp$
and that its image $\phi(\ga)$ in $G$ is regular elliptic.
Then there are two cases,
either $\phi(\ga)$ is of finite order then $\phi$ is not injective,
or $\phi(\ga)$ is of infinite order and then $\phi$ is not discrete
because the subgroup of $G$ generated by $\phi(\ga)$ is not discrete.
In particular, if a \ppp-representation $\phi$ is a discrete embedding,
then the elements $\iota_{123}$ and $\iota_{k,k-1,k,k+1}$ are not regular elliptic.

\myskip
We have shown that the trace of the element $\wA=\iota_{3231}$ is the real number 
$$\si_3=(16r_1^2r_2^2+4r_3^2-1)-16r_1r_2r_3\cos\al.$$
To check if the element~$\wA$ is regular elliptic we use the discriminant
function~$\rho$ described in section~\ref{cpxhypplane}.
As explained there, an element $A\in\SU(2,1)$ whose trace is real is regular elliptic iff $\trace A\in(-1,3)$.
The inequality $\si_k>-1$ holds for any $k=1,2,3$ and $\al\in\r$
except in the case $2r_1r_2=r_3$ and $ \al\in2\pi\cdot\z$ since
\begin{align*}
  \si_k
  &=(16r_{k-1}^2r_{k+1}^2+4r_k^2-1)-16r_1r_2r_3\cos\al\\
  &\ge16r_{k-1}^2r_{k+1}^2+4r_k^2-1-16r_1r_2r_3\\
  &=4(2r_{k-1}r_{k+1}-r_k)^2-1>-1.\\
\end{align*}
It remains to resolve the inequality $\si_3<3$ with respect to $\cos\al$.
We obtain the following proposition

\begin{prop}
\label{def-tA}
The element $\wA=\iota_{3231}$ in a \ppp-triangle group is regular elliptic iff $|t|>\tA$,
where 
$$
  \tA=\sqrt{\frac{1+\cA}{1-\cA}}
  \quad\hbox{and}\quad
  \cA=\frac{4r_1^2r_2^2+r_3^2-1}{4r_1r_2r_3}.
$$
(For $\cA>1$ resp.\ $\cA<-1$ we set $\tA=+\infty$ resp.\ $\tA=-\infty$.)
\end{prop}

\noindent
This proposition implies

\begin{cor}
\label{cor-tA}
The \ppp-triangle group representation with $|t|>\tA$ is not a discrete embedding,
because the element $\wA$ is regular elliptic.
\end{cor}

\begin{prop}
It holds $\cA\le\cinf$ and $\tA\le\tinf$.
\end{prop}

\begin{proof}
We obtain 
\begin{align*}
  4r_1r_2r_3\cdot(\cinf-\cA)
  &=r_3^2-(2r_1^2-1)(2r_2^2-1)\\
  &=\cos^2\left(\frac{\pi}{p_3}\right)-\cos\left(\frac{2\pi}{p_1}\right)\cos\left(\frac{2\pi}{p_2}\right).
\end{align*}
For $p_3\ge p_1,p_2\ge4$ this implies $\cinf\ge\cA$ because of 
$$\cos\left(\frac{\pi}{p_3}\right)\ge\cos\left(\frac{2\pi}{p_1}\right),\cos\left(\frac{2\pi}{p_2}\right)\ge0.$$
In the case~$p_1=3$ the computations are straightforward.
\end{proof}

\begin{rem}
Because of $\si_k-\si_{k+1}=4(r_k^2-r_{k+1}^2)(1-4r_{k-1}^2)$
it holds $\si_1\ge\si_2\ge\si_3$.
This means that among the elements $\iota_{k,k-1,k,k+1}$ the element
$\wA=\iota_{3231}$ is the first one to become elliptic
as $t>0$ tends to~$\tinf$ and therefore there is no need to control the
ellipticity of the elements $\iota_{1312}$ and $\iota_{2123}$.
\end{rem}

To check the ellipticity of the element $\wB=\iota_{123}$ is computationally a
much more involved task.
A tedious but straightforward computation shows that $\rho(\trace(\wB))$
is of the form
$$\rho(\trace(\wB))=\frac{f_B(t)}{(t^2+1)^3},$$
where $f_B$ is an even polynomial of degree~$6$.
Moreover, it holds $f_B(0)>0$ for $r_1,r_2,r_3>1/2$,
and the coefficient by~$t^6$ vanishes if $\cinf=1$.
The polynomial~$f_B$ can be computed explicitly,
but we shall not use these formulas, so we omit them.

\section{Triangle Groups with $r_1^2+r_2^2+r_3^2=1+2r_1r_2r_3$}

\label{spec-family}

In this section we study triangle groups corresponding
to triples $(r_1,r_2,r_3)$ with
\begin{align}
r_1^2+r_2^2+r_3^2=1+2r_1r_2r_3.\label{eq-spec-family}\tag*{$(*)$}
\end{align}

\myskip
These triangle groups seem to share a lot of properties with the ideal
triangle groups.
Possibly some of the methods used for ideal case can be used in this case too.

\myskip
The equation~\ref{eq-spec-family} is equivalent to $\cinf=1$ and $\tinf=+\infty$.
This implies that all real numbers occur as values of parameter~$t$ for some
$(\varphi_1,\varphi_2,\varphi_3)$-triangle in $\chp$.

\myskip
We first consider the case of triangle group corresponding
to a complex hyperbolic triangle with all three vertices in $\chp\cup\dd\chp$.
In this case $r_k=\cos\varphi_k$,
where $\varphi_1$, $\varphi_2$, and $\varphi_3$ are the angles of the triangle.
We assume $r_1\le r_2\le r_3$ and hence $\varphi_1\ge\varphi_2\ge\varphi_3\ge0$.
The equation~\ref{eq-spec-family} is then equivalent to
$$\varphi_1=\varphi_2+\varphi_3.$$
The signature of a non-ideal triangle group which belongs to the family~\ref{eq-spec-family}
is of the form 
$$(p_1,p_2,p_3)=(abp,a(a+b)p,b(a+b)p)$$
for some positive integer numbers $a$, $b$, and~$p$ with $a\le b$.
Some special cases are $(p,2p,2p)$-groups, $(2p,3p,6p)$-groups etc.
If we allow also ideal vertices, we obtain \ppi-groups studied by J.~Wyss-Gallifent
(compare (\ref{ppi}) in section~\ref{hist-geogr})
and among them the ideal triangle groups.
A computation shows that
$$\cA=1-\frac{\sin^2(\varphi_1+\varphi_2)}{4R}<1.$$

\myskip
In the case of an ultra-parallel complex hyperbolic
\lll-triangle group we similarly obtain
$\ell_3=\ell_1+\ell_2$ and 
$$\cA=1+\frac{\sinh^2(\ell_1-\ell_2)}{4R}>1.$$
Some special cases are $[\ell,\ell,2\ell]$-groups studied by J.~Wyss-Gallifent
(compare (\ref{ll2l}) in section~\ref{hist-geogr}), $[\ell,2\ell,3\ell]$-groups etc.

\myskip
The computation of 
$$\rho(\trace(\wB))=\frac{f_B(t)}{(t^2+1)^3}$$
for the triples $(r_1,r_2,r_3)$ satisfying the equation~\ref{eq-spec-family}
is simpler.
We obtain
$$\frac{f_B(t)}{1024R}=(1-R)t^4+(64R^3-80R^2+11R+2)t^2+(64R^3+48R^2+12R+1),$$
where $R=r_1r_2r_3$.
The polynomial $f_B$ has no roots if $R<7/8$.
For $R\ge 7/8$ the roots of the polynomial~$f_B$ can be computed explicitly
as $\pm\tBpm$, where
$$\tBpm=\sqrt{\frac{2+11R-80R^2+64R^3\pm R\sqrt{(8R-7)^3(8R+1)}}{2(R-1)}}.$$
More precisely,
\begin{enumerate}
\item for $R<7/8$ we have $f_B(t)>0$ for all~$t\in\r$,
\item for $R=7/8$ we have $f_B(t)>0$ for all $t\in\r\backslash\{\pm\tBpm\}$,
\item for $7/8<R<1$ we have $f_B(t)>0$ for $|t|<\tBm$ or $|t|>\tBp$ and
                            $f_B(t)<0$ for $\tBm<|t|<\tBp$,
\item for $R>1$ we have $f_B(t)>0$ for $|t|<\tBp$ and $f_B(t)<0$ for $|t|>\tBp$.
\end{enumerate}
In the ideal triangle group case we have $R=1$ and hence (compare~\cite{GP})
$$\rho(\tau_{123})\cdot(t^2+1)^3=1024\cdot(125-3t^2).$$

\myskip
In the ultra-parallel case the element $\wA$ remains hyperbolic for all
$t\in\r$ and the element $\wB$ goes elliptic for $t=\tBp$, hence 

\begin{prop}
Any ultra-parallel \lll-triangle group such that
$$r_1^2+r_2^2+r_3^2=1+2r_1r_2r_3,$$
\ie any ultra-parallel $[\ell_1,\ell_2,\ell_1+\ell_2]$-triangle group 
is of type B.
\end{prop}

In the non-ideal case we obtain the following result

\begin{prop}
\label{typeB}
A \ppp-triangle group such that
$$r_1^2+r_2^2+r_3^2=1+2r_1r_2r_3$$
and 
$$r_1r_2r_3\ge\frac{13+\sqrt{297}}{32}\approx 0{,}9448$$
is of type B.
For instance, the triples $(p,p,\infty)$ with $p\ge14$
and the triples $(p,2p,2p)$ with $p\ge12$ are of type B.
\end{prop}

\begin{proof}
The element~$\wA$ goes elliptic for $t=\tA$.
The element~$\wB$ goes elliptic for $t=\tBm$ since $7/8<R<1$.
We are going to prove that $\tA>\tBm$.
We have
$$\cA-R=\frac{(4r_1^2r_2^2-1)(1-r_3^2)}{4R}\ge0$$
and hence $\cA\ge R$.
This implies $\tA\ge\tAn$, where
$$\tAn=\sqrt{\frac{1+R}{1-R}}.$$
($\tAn$ is the value of $\tA$ for $r_1=r_2=\sqrt{R}$ and $r_3=1$.)
We obtain
$$f_B(\tAn)=-\frac{2048R}{1-R}\cdot(16R^2-13R-2).$$
For $R\ge(13+\sqrt{297})/32$ it holds $16R^2-13R-2>0$.
This implies $f_B(\tAn)<0$ and hence $\tAn\ge\tBm$.
The inequalities $\tA\ge\tAn>\tBm$ imply $\tA>\tBm$.
\end{proof}

\section{Applications: Arithmetic Properties of Traces}

\label{appl-disc}

Let $\Gpppn$ be the \ppp-triangle group such that $\iota_{3132}$ is a rotation by $2\pi/n$.
Then $\tau_{3132}=1+2\cos(2\pi/n)$ (compare Lemma~3.18 in~\cite{P03}).
On the other hand
$$
  \tau_{3132}
  =\si_3
  =(16r_1^2r_2^2+4r_3^2-1)-16r_1r_2r_3\cos\al
$$
and hence
$$\cos\al=(8r_1r_2r_3)^{-1}\cdot\left((8r_1^2r_2^2+2r_3^2-1)-\cos\frac{2\pi}{n}\right).$$
Here $n=\infty$ is allowed and means $\tau_{3132}=3$.

\begin{prop}
\label{uuu-traces}
Let $\tau$ be the trace of the element in $\Gpppn$, then
$$
  2\Re(\tau),~|\tau|^2\in
  \z
  \left[
        2\cos\left(\frac{2\pi}{p_1}\right),
        2\cos\left(\frac{2\pi}{p_2}\right),
        2\cos\left(\frac{2\pi}{p_3}\right),
        2\cos\left(\frac{2\pi}{n}\right)
  \right].
$$
\end{prop}

\begin{proof}
Because of $2\cos(2\pi/p_k)=4r_k^2-2$
we have
$$
  \z\left[
          2\cos(2\pi/p_1),
          2\cos(2\pi/p_2),
          2\cos(2\pi/p_3),
          2\cos(2\pi/n)
    \right]
  =R\left[2\cos(2\pi/n)\right],
$$
where $R=\z[4r_1^2,4r_2^2,4r_3^2]$.
For the trace $\tau$ of an element of $\Gpppn$ it holds
according to the combinatorial trace formula $\tau=(-1)^{|a|}\cdot(2+\xi)$,
where $|a|$ is the length of $a$ and
$$\xi=\sum_{w\in\z} q_w\cdot e^{i\al w}.$$
We have $2\Re(\tau)=(-1)^{|a|}\cdot(4+2\Re(\xi))$ and $|\tau|^2=4+4\Re(\xi)+|\xi|^2$,
hence is is sufficient to prove that $2\Re(\xi)$ and $|\xi|^2$
belong to the ring $R[2\cos(2\pi/n)]$.
We have
$$2\Re(\xi)=\sum_{w\in\z} q_w\cdot(2\cos(\al w))$$
and
$$
  |\xi|^2
  =\xi\cdot\bar\xi
  =\sum_{w\in\z} q_w^2+\sum_{{w,w'\in\z\atop w\ne w'}} q_w\cdot q_{w'}\cdot(2\cos(\al(w-w'))).
$$
But $2\cos(\al w)$ is a polynomial of degree~$|w|$ in $2\cos\al$
with integer coefficients,
even for even $|w|$ and odd for odd $|w|$
(Chebyshev polynomials, compare for example~\cite{S95}).
This implies
$$2\cos(\al w)\in(8r_1r_2r_3)^{-|w|}\cdot R[2\cos(2\pi/n)]$$
because of
$$2\cos\al\in(8r_1r_2r_3)^{-1}\cdot R[2\cos(2\pi/n)]$$
and hence
\begin{align*}
  (2\cos\al)^{|w|-2j}
  &\in(8r_1r_2r_3)^{-|w|}\cdot(4r_1^2\cdot 4r_2^2\cdot 4r_3^2)^j\cdot R[2\cos(2\pi/n)]\\
  &\subset(8r_1r_2r_3)^{-|w|}\cdot R[2\cos(2\pi/n)].\\
\end{align*}
On the other hand the coefficients $q_w$ satisfy $q_w\in(8r_1r_2r_3)^{|w|}\cdot R$ by lemma~\ref{qw-in-ring}.
This finishes the proof.
\end{proof}

\begin{cor}
Let $\tau$ be the trace of an element in $\Gpppn$ with
$$\{p_1,p_2,p_3,n\}\subset\{3,4,6,\infty,q\},$$
where $q$ is any natural number, 
then $2\Re(\tau)$ and $|\tau|^2$ both belong to the ring $\z[2\cos(2\pi/q)]$.
\end{cor}

\begin{proof}
For $a\in\{3,4,6,\infty\}$ we have $2\cos(2\pi/a)\in\z$.
\end{proof}

\section{History and Geography}

\label{hist-geogr}

In this section we describe some of the results in the area of
complex hyperbolic triangle groups and specify there position in the space of
triples $(r_1,r_2,r_3)$.

\begin{enumerate}[(1)]
\item ideal triangle case: $r_1=r_2=r_3=1$, ($p_k=\infty$, $\varphi_k=0$), type B.

\myskip\noindent
Ideal triangle group representations were first studied in~\cite{GP}
by William Goldman and John Parker.
They proved that an ideal triangle group representation is a discrete embedding
for $\cos\al<\frac{17}{18}$ and is not a discrete embedding for
$\cos\al>\frac{61}{64}$, and they conjectured that an ideal triangle group
representation is a discrete embedding for $\cos\al\le\frac{61}{64}$.
This conjecture was confirmed by Richard E.~Schwartz in~\cite{Sch2001:dented}.
He proved that an ideal triangle group representation is a discrete embedding
if and only if $\cos\al\le\frac{61}{64}$.
In~\cite{Sch2001:degen} he also used the so called last ideal triangle group, \ie the group with
$\iota_{123}$ parabolic, \ie with $\cos\al=\frac{61}{64}$,
to construct the first example of a complete hyperbolic $3$-manifold with a spherical CR-structure.
Hanna Sandler suggested in~\cite{S95} to study formulas for the traces of
elements in an ideal triangle group in order to prove the Goldman-Parker conjecture.

\myskip
\item \uuu-groups: $r_1=r_2=r_3=2^{-1/2}$, ($p_k=4$, $\varphi_k=\pi/4$), type A.

\myskip\noindent
Richard E.~Schwartz studied in~\cite{Sch2003} the groups $\Guuun$ and used the
group $G(4,4,4;7)$, \ie the \uuu-group with $\cos\al=2^{-3/2}\cdot(2-\cos(2\pi/7))$
to construct the first example of a compact hyperbolic $3$-manifold with a spherical CR-structure.

\myskip
\item \uui-groups: $r_1=r_2=2^{-1/2}$, $r_3=1$, type A.

\myskip\noindent
Justin Wyss-Gallifent studied \uui-groups in~\cite{W}.
He proved that a \uui-group is not discrete for $\cos\al>1/2$
and that the groups $\Guuin$ are discrete for $n=3,4,\infty$.

\myskip
\item \ppi-groups: $r_1=r_2$, $r_3=1$.

\label{ppi}

\myskip\noindent
The \ppi-groups were studied also by Justin Wyss-Gallifent in~\cite{W}.
He proved that a triple \ppi\ is of type~A for $p\le13$ and of type~B for $p\ge14$.

\myskip
\item ultra-parallel $[\ell,\ell,2\ell]$-groups: $r_1=r_2=r$, $r_3=2r^2-1$.

\label{ll2l}

\myskip\noindent
Justin Wyss-Gallifent studied in~\cite{W} ultra-parallel $[\ell,\ell,2\ell]$-groups
and obtained partial results similar to the results in~\cite{GP} on ideal
triangle groups.

\myskip
\item Mostow's non-arithmetic complex hyperbolic lattices:
$\mu_1=\mu_2=\mu_3=e^{\frac{2\pi i}{p}}$, $p\in\{3,4,5\}$,
$r_1=r_2=r_3=\frac{1}{2\sin\frac{\pi}{p}}$,
$\al=\frac{2\pi}{\rho}+\frac{\pi}{p}-\frac{\pi}{2}$, $\rho\in\z$.

\myskip\noindent
G.~D.~Mostow studied in~\cite{M} some triangle groups generated by three $\mu$-reflections
of the same order $p\in\{3,4,5\}$ in the sides of a equiangular triangle with
angles related to the order $p$ and with special values of the angular invariant $\al$
and obtained the first examples of non-arithmetic complex hyperbolic lattices.
\end{enumerate}


\bibliographystyle{amsalpha}
\bibliography{formula2}

\end{document}